\documentclass[graybox]{svmult}

\usepackage{mathptmx}
\usepackage{helvet}  
\usepackage{courier} 
\usepackage{type1cm} 
                     
\usepackage{makeidx} 
\usepackage{graphicx}
                     
\graphicspath{{Figures/}}
\usepackage{multicol} 
\usepackage[bottom]{footmisc}
\usepackage{amsmath}

\makeindex 

\begin{document}
	
\title*{Parametric identification of the dynamics\\ of inter-sectoral 
balance: modelling and forecasting\thanks{This is a preprint 
of a paper accepted for publication 29-March-2019 
as a book chapter in ``Advances in Intelligent Systems and Computing'' 
(\url{https://www.springer.com/series/11156}), Springer.}}

\titlerunning{Parametric identification of the dynamics of inter-sectoral balance}

\author{Olena Kostylenko, Helena Sofia Rodrigues and Delfim F. M. Torres}

\authorrunning{O. Kostylenko, H. S. Rodrigues and D. F. M. Torres}

\institute{Olena Kostylenko (corresponding author) 
\at Center for Research and Development in Mathematics and Applications (CIDMA),\\ 
Department of Mathematics, University of Aveiro, 3810-193 Aveiro, Portugal\\ 
\email{o.kostylenko@ua.pt}
\and 
Helena Sofia Rodrigues 
\at Center for Research and Development in Mathematics and Applications (CIDMA),\\ 
Department of Mathematics, University of Aveiro, 3810-193 Aveiro, Portugal;\\ 
School of Business Studies, Polytechnic Institute of Viana do Castelo, 4930-678 Valen\c{c}a, Portugal 
\email{sofiarodrigues@esce.ipvc.pt}
\and 
Delfim F. M. Torres 
\at Center for Research and Development in Mathematics and Applications (CIDMA),\\ 
Department of Mathematics, University of Aveiro, 3810-193 Aveiro, Portugal\\ 
\email{delfim@ua.pt}}

\maketitle
	

\abstract{This work is devoted to modelling and identification 
of the dynamics of the inter-sectoral balance of a macroeconomic system. 
An approach to the problem of specification and identification 
of a weakly formalized dynamical system is developed. A matching procedure
for parameters of a linear stationary Cauchy problem with a decomposition of its upshot 
trend and a periodic component, is proposed. Moreover, an approach for detection 
of significant harmonic waves, which are inherent to real 
macroeconomic dynamical systems, is developed.}

\bigskip


\medskip


	
\section{Introduction}
\label{sec:1}

For any dynamical object (technical, economic, environmental, etc.), 
the most critical is the problem of limited resources and their optimal use. 
Thereby, there is a need to build mathematical models that adequately describe 
the existing trends and provide high-precision predictive properties of the 
dynamical systems. Therefore, it is necessary to evaluate the quality 
of mathematical models of dynamical processes through the prism 
of imitation and predictive properties. 

The meaning of mathematical modelling, as a research method, is determined 
by the fact that a model is a conceptual tool focused on the analysis 
and forecasting of dynamical processes using differential 
or differential-algebraic equations. However, the parameters of these 
equations are usually unknown in advance. Therefore, in practice, any 
direct problem (imitation, prediction, and optimization) is always preceded 
by an inverse problem (model specification and identification 
of parameters and variables included in it).

To construct dynamical models in economics, it is necessary to 
formulate principles, based on economic theory, and consider
equations that adequately determine the evolution of the investigated process. 
The functioning of a market economy, as well as any economic system, 
is not uniform and uninterrupted. An important tool for its prediction 
is the dynamical model of inter-industry equilibrium balance developed by Leontief 
(the so called ``Leontief's input-output model''), 
which became the basis of mathematical economics \cite{[Leontief]}. 
Leontief's input-output model is described by a system of linear 
differential equations in which the input is the non-productive consumption of sectors, 
and the output is the issue of these sectors. In practice, the matrices of the system 
of differential equations are unknown in advance. Therefore, there is a need to evaluate them 
based on statistical information of the macroeconomic system under study, 
for a certain period of time. The problem of parametric identification of macroeconomic 
dynamics is discussed in the papers \cite{[Ramsay],[NazarenkoFil],Olena01,Olena02}.

One feature of the dynamical Leontief model is that it can be considered 
as taking into account the cyclical nature of the macroeconomic processes, 
typical for developed countries of Western Europe, USA, Canada, Japan, 
and others \cite{[Diebolt],[Korotaev]}. According to statistical analysis, 
their economic growth (upturn phase) alternates with the processes of 
stagnation and decline in production volumes (downturn phase), i.e., 
a decline in all the economic activity. After a decline phase, again, 
the economic cycle continues, similarly with ups and downs.
Such periodic fluctuations indicate a cyclical nature of the economic development.

The analysis of the input-output model is critical for developing 
government regulatory programs, that are an integral part of the 
economic and environmental policy for any developed market system. 
Such an analysis is a useful tool for studying the strategic directions 
of economic development, the expansion of dynamical inter-industry relations, 
and the interaction between economic sectors and the natural environment 
\cite{[Fied],[Dobos]}. 

This paper is organized as follows. Section~\ref{sec:2} is devoted to a 
detailed economic and mathematical description of the Leontief macroeconomic model. 
Section~\ref{sec:3} focuses on the problem of specification and identification 
of this model and Section~\ref{sec:4} to the testing of the Leontief model with 
time series data of a real macroeconomic dynamics. 
In Section~\ref{sec:5}, the main conclusions are carried out.


\section{Problem statement}
\label{sec:2}

The static model of inter-industry balance is obtained when one equates 
the net production of branches to the final demand for the products 
of these industries:
\begin{equation}
\label{eq:01}
\textbf{x} - A\textbf{x} = \textbf{y}\;,
\end{equation}
where ${\bf x}=(x_{1},x_{2},\ldots,x_{n})^{'}$ is the column vector 
of gross output by industries; ${\bf y}=(y_{1},y_{2},\ldots,y_{n})^{'}$ 
is the column vector of final demand; 
and $A$ is the matrix of direct costs. The vector of final demand ${\bf y}$ 
can be represented as a sum of two vectors: investments, 
${\bf y}_{1}(t)=B\dot{\bf x}(t)$, and consumer products, ${\bf y}_{2}(t)={\bf u}(t)$. 
Then, the model of dynamical inter-industry balance is
\begin{eqnarray}
\label{eq:02}
{\bf x}(t)=A{\bf x}(t) +B\dot{\bf x}(t) + {\bf u}(t)\;, \ \
\dot{\bf x}(t)=\frac{d{\bf x}(t)}{dt}\;,
\end{eqnarray}
where ${\bf u}(t)$ is the vector-column of final (non-productive) consumption 
and $B$ is the matrix of capital coefficients.
A discrete-time dynamic balance model can be reduced 
to a model with continuous time as follows:
\begin{eqnarray}
B\dot{\bf x}(t) = (E-A){\bf x}(t)-{\bf u}(t)\;.
\label{eq:03}
\end{eqnarray}
If the inverse matrix $B^{-1}$ exists, 
then the Leontief's model can be rewritten as a linear multiply-connected system, 
with the final consumption vector ${\bf u}(t)$ as the input 
and with the gross output vector ${\bf x}(t)$ as the output:  
\begin{equation}
\label{eq:04}
\dot{\bf x}(t) = B^{-1}(E-A){\bf x}(t) - B^{-1}{\bf u}(t)\;, \ \
t\in [t_0,t_f]\;.
\end{equation}
The differential equation \eqref{eq:04} must be supplemented 
with the boundary condition
\begin{equation}
\label{eq:05}
{\bf x}(t_*) = {\bf x}_*\;,
\end{equation}
where $t_{*}$ is the point of the segment $[t_{0},t_{f}]$ 
at which the boundary state is specified. 
To determine the unknown parameters
$A$, $B$, and ${\bf x}_*$ of the Cauchy problem 
\eqref{eq:04}--\eqref{eq:05}, in accordance with the methodology 
proposed in \cite{[Nazarenko]}, we divide the interval 
$[t_{0},t_{f}]$ into two intervals: the identification period
$[t_{0},t_{*})$ and the forecasting period $[t_{*},t_{f}]$. 
We assume that the base period $[t_{0},t_{*})$ has statistical 
information {${\bf x}_t$} and {${\bf y}_t$} 
regarding industry's outputs and consumption at its discrete-times 
$t = 1, 2,\ldots, N$, and the model will be adjusted for the vector 
of unknown coefficients ${\bf \theta}$. For the interval $[t_{*},t_{f}]$, 
we assume that $t_{f}-t_{*} \ll N$. Through this approach, the boundary 
condition \eqref{eq:05} is satisfied in the first integer point of the forecast period.
Then, if the model is stationary, and due to the inertia of the dynamical system, 
the vector ${\bf \theta}$ can be translated for the forecast period. 
The stationariness of the model is characterized by the high quality 
of approximation and forecasting and robustness \cite{[Green]}.
However, for modelling and prediction of linear stationary models, 
the identification period $N$ should be large enough to 
stabilize the interrelations between the elements of the system.


\section{Algorithm}
\label{sec:3}

First of all, for convenience, we make the data dimensionless and then, 
we proceed by solving the Cauchy problem \eqref{eq:04}--\eqref{eq:05}.
To do this, we need to specify the vector of phase coordinates ${\bf x}(t)$ 
and estimate the unknown parameters that will appear in the specification 
process \cite{[Nazarenko]}. The vector ${\bf x}(t)$ is searched 
by decomposing the trajectories of the phase coordinates motion 
into components \cite{[Nazarenko]}. 
If these trajectories are defined, 
then the vector ${\bf u}(t)$ can be found.
It will be supplied to the input of the dynamical system 
and it can be used to solve the problem of specifying the sectors' 
output and manage their movement. The vector ${\bf u}$ will be 
considered as a control vector and it can be found 
using the inverse relationship
\begin{equation}
\label{eq:06}
{\bf u}(t) = P{\bf x}(t)-B\dot{\bf x}(t)\;, \ \ P=E-A\;.
\end{equation}
At each time step, the vector $\bf{u}$ is a function 
of the phase coordinates and their derivatives. 
The regulator implements the critical idea of control theory: 
the principle of inverse connection, being used for identifying 
the model \eqref{eq:06}. In our research, the regulator consists of two devices.
At any time $t$, they form a total value of phase coordinates (issues sectors) 
and controls (non-productive consumption sectors):
\begin{eqnarray}
\label{eq:07}
{\bf x}(t)=\sum \limits_{m=1}^{n}{\bf x}_m(t)\;, \ \
{\bf u}(t)=\sum \limits_{m=1}^{n}{\bf u}_m(t)\;.
\end{eqnarray}
If the modelling trajectories ${\bf x}(t)$ and ${\bf u}(t)$, where
\begin{math}
{\bf x}(t)=(x_{1}(t),x_{2}(t),\ldots,x_{n}(t))^{'}
\end{math}
and 
\begin{math}
{\bf u}(t)=(u_{1}(t),u_{2}(t),\ldots,u_{n}(t))^{'},
\end{math} 
are adjusted to high imitation and prediction properties, 
then the total trajectories \eqref{eq:07} should have the same properties.
Following \cite{[Nazarenko]}, we assume that the trajectory of the 
dynamical system is represented by an additive combination of its components:
\begin{equation}
\label{eq:08}
{\bf x}_t = {\bf x}_*(t)+\varepsilon_t\;.
\end{equation}
The development tendency is characterized by a linear trend,
\begin{eqnarray}
\label{eq:9}
{\bf x}_*(t) = \bar{\bf x}+{\bf b}(t-\bar{t})\;, \ \
\bar{\bf x}_* = \bar{\bf x}+{\bf b}(t-\bar{t})=\bar{\bf x}\;, 
\end{eqnarray}
while the fluctuation process is described by a linear combination 
of harmonics with some frequencies over the time interval $[1, N]$:
\begin{equation}
\label{eq:10}
{\bf \varepsilon}_t  = \sum \limits_{k}({\bf a}_k 
\cos \omega_{k}t + {\bf b}_k \sin \omega_{k}t) +\nu_t\;, 
\end{equation}
where $\omega_{k}$ is the frequency of the harmonic $k$; 
${\bf a}_{k},{\bf b}_{k}$, $k=1,\ldots,n-1$, are the 
unknown coefficients of the decomposition in the truncated Fourier series; 
and ${\bf \nu}_{t}$ is the vector of random residuals. We decompose the process 
of random oscillations as a truncated Fourier series in order to reflect 
the imitation and predict the real properties of the oscillations.
Since the trend and periodic components correlate with each other, 
then the phase trajectories can be found after evaluating the regression model: 
\begin{equation}
\label{eq:11}
{\bf x}_t -{\bf x} = {\bf b}(t-\bar{t})+ \sum \limits_{k=1}^{n-1}({\bf a}_k 
\cos \omega_{k}t + {\bf b}_k \sin \omega_{k}t) +{\bf \nu}_t\;, 
\quad t=1,\ldots, N\;.
\end{equation}
It is assumed that in regression models the average value of residuals is zero. 
Therefore, it is necessary to choose the frequency fluctuation, so that the 
average harmonics values would be equal to zero. 
Harmonic fluctuations are adjusted in the spectrum of frequencies:
\begin{eqnarray}
\label{eq:12}
\left\{
\begin{array}{ccc}
\sin \omega_{k}t & = & 0 \\
\cos \omega_{k}t & = & 0 \\
\end{array}
\right.
\Rightarrow  
\omega_{k} = \frac{2\pi}{N}k, \ \ k=1,2, \ldots
\end{eqnarray}
Determination of frequencies from the spectrum given in \eqref{eq:12}, 
as well as the period $T$ of oscillations of this system, can be done 
using the first regulator device \eqref{eq:06}, which calculates 
the output of the entire macroeconomic system. If the fluctuation 
frequencies belong to the spectrum \eqref{eq:12}, where the variable 
$N$ is the sample size, then, from the mathematical point of view, 
we need the fluctuation period $T$ instead of the sample size. 
So, we assume that $N=T$. The optimal value $N$, as well as the 
fluctuation period, can be identified using the back-extrapolation method.
Minimization of the residual sum of squares 
\begin{equation}
\label{eq:13}
S=\sum_{t=1}^{N} \nu_{t}^2 = \sum_{t=1}^{N}(x_t-\bar{x}-b(t-\bar{t})
-\sum \limits_{k}^{}(a_k \cos \omega_{k}t + b_k \sin \omega_{k}t))^2\; 
\end{equation}
provides the following Ordinary Least Squares (OLS) estimates:
\begin{equation}
\label{eq:14}
\begin{gathered}
\hat{b}= \frac{\frac{2}{N} \sum \limits_{t=1}^{N} {(x_t - \bar{x})(t
+ \sum \limits_{k}(\cot \frac{\omega_{k}}{2} \sin \omega_{k}t - \cos \omega_{k}t))}}
{\frac{N^2 -1}{6} -\sum \limits_{k} 
\frac{1}{\sin^{2} \frac{\omega_{k}}{2} }} \;,  \\
\varepsilon_t  = x_t -\bar{x}-\hat{b}(t-\bar{t})\;, \ \
\hat{a_k}= \frac{2}{N} \sum \limits_{t=1}^{N} x_t \cos \omega_k t - \hat{b}\;,\\
\hat{b_k}= \frac{2}{N} \sum \limits_{t=1}^{N} x_t \sin \omega_k t 
+ \hat{b} \cot \frac{\omega_k}{2}\;, \ \
\hat{c_k}= \sqrt{\hat{a_k}^2 + \hat{b_k}^2} \;.
\end{gathered}
\end{equation}
According to \eqref{eq:14}, the trend correlates with harmonics and, vice versa, 
harmonic waves depend on the trend around which they fluctuate. 
Harmonics of the regression model of fluctuations, which have different 
spectra \eqref{eq:12}, do not correlate with each other.
If the estimates $\hat a_{k}$ and $\hat b_{k}$, $k=1,2,\ldots$, are optimal
\eqref{eq:14}, then the minimum sum of squares of deviations is equal to
\begin{equation}
\label{eq:15}
\min S = \sum \limits_{t=1}^{N} \hat{\nu}^2 
= \sum \limits_{t=1}^{N} \varepsilon_t ^2 
- \frac{N}{2} \sum \limits_{k} \hat{c_k}^2 \;.
\end{equation}
The significance of the OLS-estimates \eqref{eq:14} 
was tested using Student's t-test.
The optimal dimension $n$ of the phase space is determined 
when the $n-1$ significant harmonics are found.
The component specification  ${\bf x}$ is executed with a specified value $n$.
If the phase coordinates are chosen, then we assume that the inherent harmonic 
fluctuations are tuned to the frequency \eqref{eq:12}. The number of significant 
harmonics in models of different phase coordinates may differ, due to the fact 
that each subset $x_1, x_2, \ldots, x_n$ of the set ${\bf x}$ has its own 
specificity of functioning. If the phase coordinate responds quickly to 
qualitative changes in a given dynamical system, then this coordinate will 
have the maximum number, that is, $n-1$ harmonics. The minimum number 
of harmonics will be in the decomposition of those phase coordinates 
that are weakly responsive to changes in other subsets of the system.
If the insignificant OLS-estimates of the decomposition coefficients are discarded, 
then the model of the phase coordinate fluctuations around the corresponding trends 
is back-extrapolated for the period $t \leq 0$ and checked if matches the statistical 
value $\varepsilon_{t}$, $t=0,-1,\ldots$ If we are satisfied with the test, 
then we model the trajectories of the phase coordinate movement as
\begin{equation}
\label{eq:16}
{\bf x}(t) = \bar{\bf x}+{\bf b}(t-\bar{t})+\sum \limits_{k}(\hat{\bf a}_k 
\cos \omega_{k}t + \hat{\bf b}_k \sin \omega_{k}t)\;.
\end{equation}
The approximation properties of the obtained model curves are 
described by using the coefficients of determination $R^2$.
The coefficients of determination linear trends of 
the phase coordinates are calculated by the formula
\begin{equation}
\label{eq:17}
{\bf R}_{tr}^{2} = \frac{\hat{\bf b}^2 
\frac{(N^2 -1)}{12}N}{\sum \limits_{t=1}^{N} ({\bf x}_t - \bar {\bf x})^2}\;.
\end{equation}
Taking into account the periodic component in the
model motion paths \eqref{eq:16}, we get 
\begin{equation}
\label{eq:18}
{\bf R}^{2} =1- \frac{\sum \limits_{t=1}^{N} {\bf \nu}_t ^2}{\sum 
\limits_{t=1}^{N} ({\bf x}_t - \bar {\bf x})^2}\;, \ \ \
{\bf \nu}_t ^2=\sum_{t=1}^{N} {\bf \varepsilon}_t ^2 - \frac{N}{2} \sum_{k} {\bf c}_k ^2 \;.
\end{equation}
For the fluctuation model, one gets
\begin{equation}
\label{eq:19}
{\bf R}_{fl}^{2} = \frac{ \frac{N}{2} \sum \limits_k 
\hat{\bf c}_k ^2}{\sum \limits_{t=1}^{N} {\bf \varepsilon}_t^2}.
\end{equation}
Particles of harmonic dispersions, in the general dispersion fluctuations 
of each sector $m=1, \ldots, n$, are calculated by using respective 
coefficients of determination. The proportion of dispersion of the 
$k$th harmonic, in the total dispersion of the fluctuation 
of the phase coordinates, is
\begin{equation}
\label{eq:20}
{\bf R}_{k}^{2} = 
\frac{ \frac{N}{2} \hat{\bf c}_k ^2}{\sum \limits_{t=1}^{N} {\bf \varepsilon}_t^2} 
=\frac{\hat{\bf c}_k ^2}{\sum \limits_{k}\hat{\bf c}_k ^2} {\bf R}_{fl}^{2} \;.
\end{equation}
The specification of the control vector ${\bf u}$ is carried out 
at the given phase vector ${\bf x}(t)$ and its derivative $\dot{\bf x}(t)$.
To identify the control vector ${\bf u}(t)$ taking into account \eqref{eq:06}, 
we construct the following regression model: 
\begin{equation}
\label{eq:21}
{\bf u}(t)-\bar{\bf u} = P\left({\bf x}(t)-\bar{\bf x}\right) 
- B\left(\dot{\bf x}-\hat{\bf b}\right) + {\bf r}_t\;,
\end{equation}
where ${\bf r}_t$ is the vector of random perturbations.
The adequacy of the curves is checked using the coefficients of determination, 
as well as the second controller device, which models the trajectory of 
total non-production consumption, according to the second balance equation \eqref{eq:07}.


\section{Implementation}
\label{sec:4}

The approbation of the constructed algorithm was carried out for 
real macroeconomic dynamics. The research was done using 
the statistical data of France \cite{[France]}, 
which includes information from 1949 to 2017.
It was found that the period of fluctuation for France is fifty years $(T=50)$.
Consequently, $N=50$, while 1966--2015 is the identification period 
and 2016--2017 the forecast period. The next and main step 
is to establish the significant harmonics. For this, we have used 
the Student t-test with a significance level of $\alpha=0.005$. 
If the $k$-harmonic appears insignificant by Student's t-test, 
then all insignificant harmonics are sequentially turned over, 
and we need to remove them from the regression.  
Then, coming back to the previous step, we start recalculations, 
as long as all the harmonics become significant.
When we find four significant harmonics, Kondratieff's wave $(k=1)$, 
Kuznets' wave $(k=3)$, Juglar's wave $(k=6)$, and the wave that equals 
half the period of the Kondratieff wave $(k=2)$, then we conclude that the number 
of all sectors is five, that is, $n=4$ harmonics plus $1$ trend \cite{[Korotaev]}.
Thus, the economy must be divided into five sectors.
Empirically, we set the optimum division of the economy into sectors:
Industry and Agriculture (Sector~1); 
Construction and Transport (Sector~2); 
Finance and Real estate (Sector~3);
Communication and Science (Sector~4); 
and Service Industries (Sector~5); 
In each sector, the four harmonics are present.
Parametric identification of the regression model of the sectors outputs
gives the following coefficients of trends determination $R^2$ 
(Table~\ref{tab:1}), around which fluctuation occurs.
\begin{table}
\begin{center}
\caption{The coefficients of trends determination.}
\label{tab:1} 
\begin{tabular}{p{1.4cm}p{1.4cm}p{1.4cm}p{1.4cm}p{1.4cm}p{1.4cm}p{1.4cm}}
\hline\noalign{\smallskip}
 Sector & 1 & 2 & 3 & 4 & 5 & $\sum$ \\
\noalign{\smallskip}\svhline\noalign{\smallskip}
$R^2$ & 0.8076 & 0.6871 & 0.7773 & 0.7233 & 0.8151 & 0.7883\\
\noalign{\smallskip}\hline\noalign{\smallskip}
\end{tabular}
\end{center}
\end{table}
Results in Table~\ref{tab:1} show that, for the economy of France, 
the fluctuation of outputs around the corresponding trend is perceptible.
The harmonic waves interact with the trend, complicating their analysis. 
But if we consider the pure oscillatory process, then the harmonics 
of Fourier's series become uncorrelated, simplifying the analysis 
of the influence of individual harmonics on the general oscillation process.
Particles of harmonic dispersions in the general dispersion fluctuations of each sector 
$m=1,\ldots, n$ are calculated using the respective coefficients of determination 
\eqref{eq:20}, the values of which are given in Table~\ref{tab:2}, where the ``---'' 
means that the harmonic is insignificant in the sector according to Student's t-test.
\begin{table}
\begin{center}
\caption{Contribution of harmonics into the oscillatory process.}
\label{tab:2} 
\begin{tabular}{p{1.5cm}p{1.5cm}p{1.5cm}p{1.5cm}p{1.5cm}p{1.5cm}}
\hline\noalign{\smallskip}
Sector & $k=1$ & $k=2$ & $k=3$ & $k=6$ & $\sum$ \\
\noalign{\smallskip}\svhline\noalign{\smallskip}
1 & 0,8202 & 0,0254 & 0,0908 & --- & 0,9364 \\
2 & 0,3719 & 0,4622 & 0,0982 & 0,0173 & 0,9496 \\
3 & 0,6460 & 0,2916 & 0,0208 & --- & 0,9583 \\ 
4 & 0,8855 & 0,0884 & 0,0073 & --- & 0,9812 \\
5 & 0,7354 & 0,2521 & --- & --- &	0,9875 \\
$\sum$ & 0,7429 & 0,2186 & 0,0094 & 0,0089 & 0,9798 \\
\noalign{\smallskip}\hline\noalign{\smallskip}
\end{tabular}
\end{center}
\end{table}
The analysis of the share of harmonics in the general dispersions 
fluctuations of each sector, calculated by appropriate 
coefficients of determination, shows that:
\begin{itemize}
\item The Kondratieff wave $k=1$ (long wave) 
substantially influences the $1$st and $4$th sectors.

\item The wave with period of 25 years ($k=2$) 
prevails in the $2$nd and $3$rd sectors.

\item The Kuznec wave $k=3$ (rhythms with period 15--20 years), 
especially manifests itself in the $1$st and $3$rd sectors.

\item The contribution of the Juglar wave ($k=6$) is less significant
in comparison with other waves, but it makes significant adjustments 
to the gross output function for the Sectors~2.

\item The total contribution of harmonics, into the dispersion 
of fluctuations of the sectors, ranges from $93.64\%$ ($1$st sector) 
to $98.75\%$ ($5$th sector).
\end{itemize}
Therefore, regression models of fluctuations have the qualitative 
approximation properties and we can expect a significant 
contribution into the dispersions of issues.
The values of the coefficients of determination of modelling 
trajectories of issues are given in Table~\ref{tab:3}.
\begin{table}
\begin{center}
\caption{The quality of modelling trajectories of issues.}
\label{tab:3} 
\begin{tabular}{p{1.4cm}p{1.4cm}p{1.4cm}p{1.4cm}p{1.4cm}p{1.4cm}p{1.4cm}}
\hline\noalign{\smallskip}
Sector & 1 & 2 & 3 & 4 & 5 & $\sum$ \\
\noalign{\smallskip}\svhline\noalign{\smallskip}
$R^2$ & 0,9983 & 0,9972 & 0,9987 & 0,9991 &	0,9998 & 0,9995 \\
\noalign{\smallskip}\hline\noalign{\smallskip}
\end{tabular}
\end{center}
\end{table}
Comparing the results shown in Table~\ref{tab:1} and Table~\ref{tab:3}, 
we can make a conclusion about the significant influence of harmonic waves 
on trajectories of sectors issues (the quality of modelling trajectories exceeds $99\%$).

\begin{figure}[t]
\sidecaption[t]
\label{f1a}\includegraphics[width=0.5\textwidth]{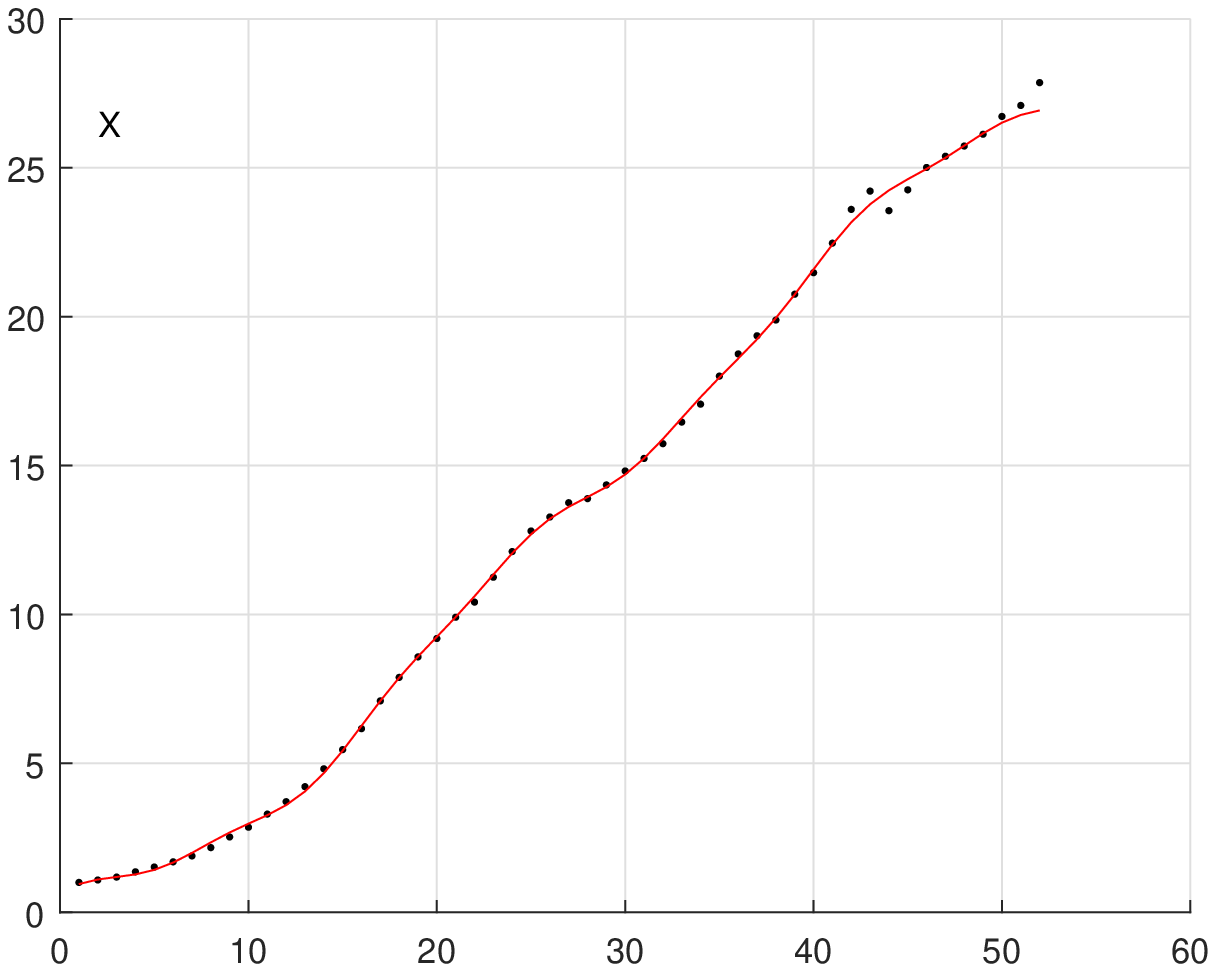}
\label{f1b}\includegraphics[width=0.5\textwidth]{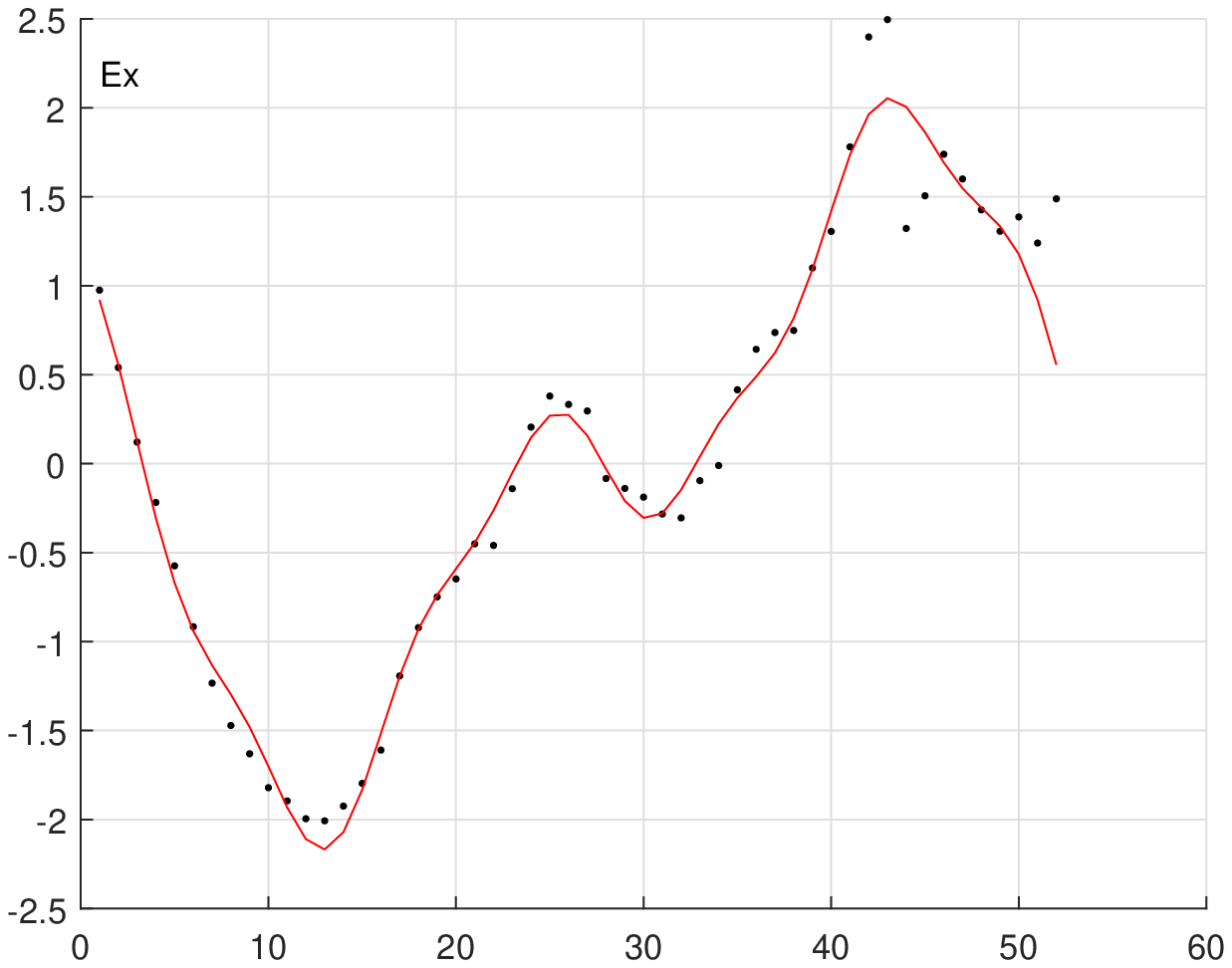}
\caption{Modelling curves of gross outputs (left) 
and fluctuation of gross output (right) of all economy.}
\label{fig:1}   
\end{figure}

In Figure~\ref{fig:1}, the graphics of modelling curves of gross issues 
and the trajectories corresponding to the fluctuations 
of the macroeconomic system are shown.
The graphics show points that give statistical information 
and the solid line is the trajectory. Comparison of the predicted 
values with the real data (the last two points, which correspond 
to years $2016$ and $2017$), testifies that the model 
has high accuracy predictive properties. 

\begin{figure}[t]
\sidecaption[t]
\label{f2a}\includegraphics[width=0.5\textwidth]{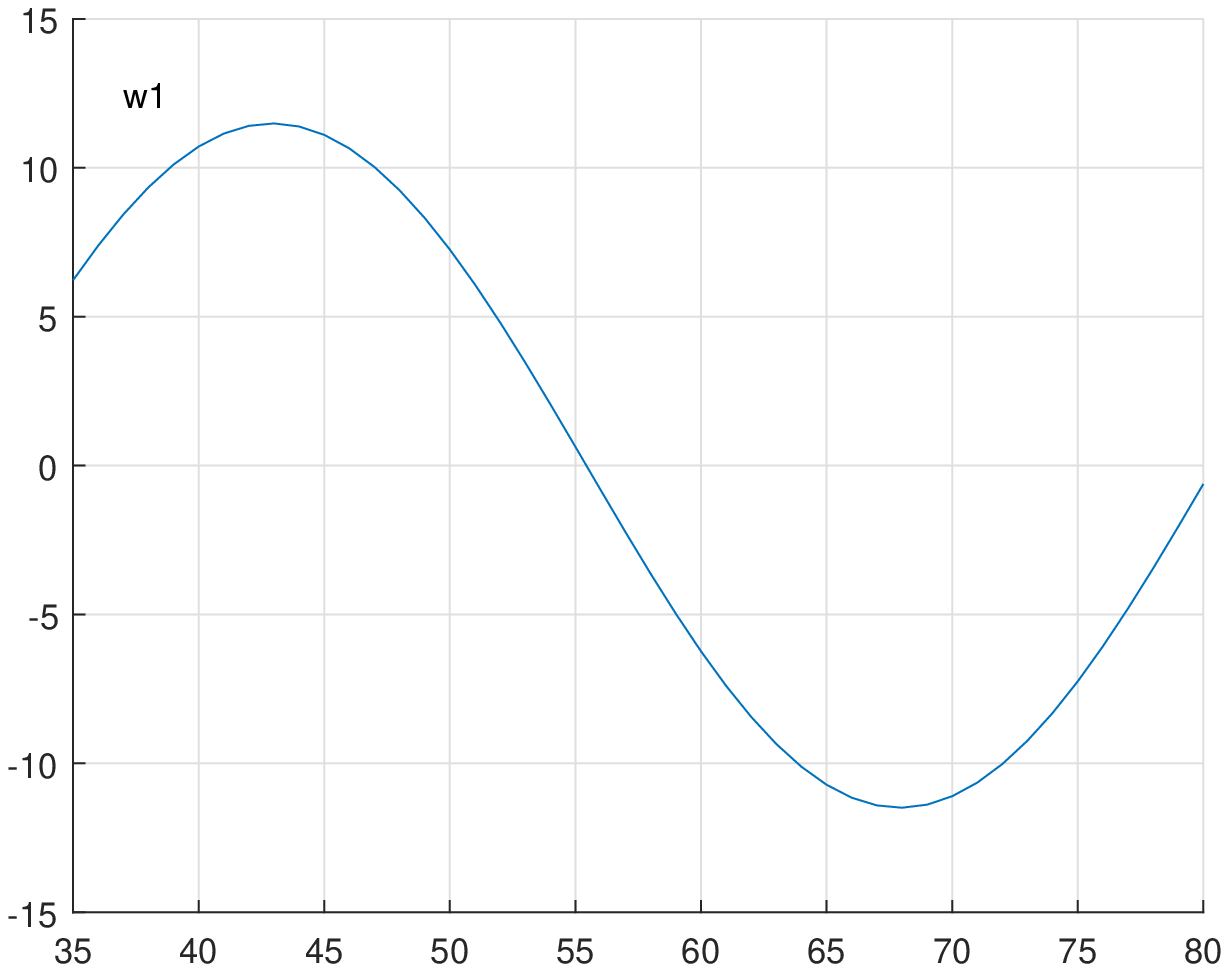}
\label{f2b}\includegraphics[width=0.5\textwidth]{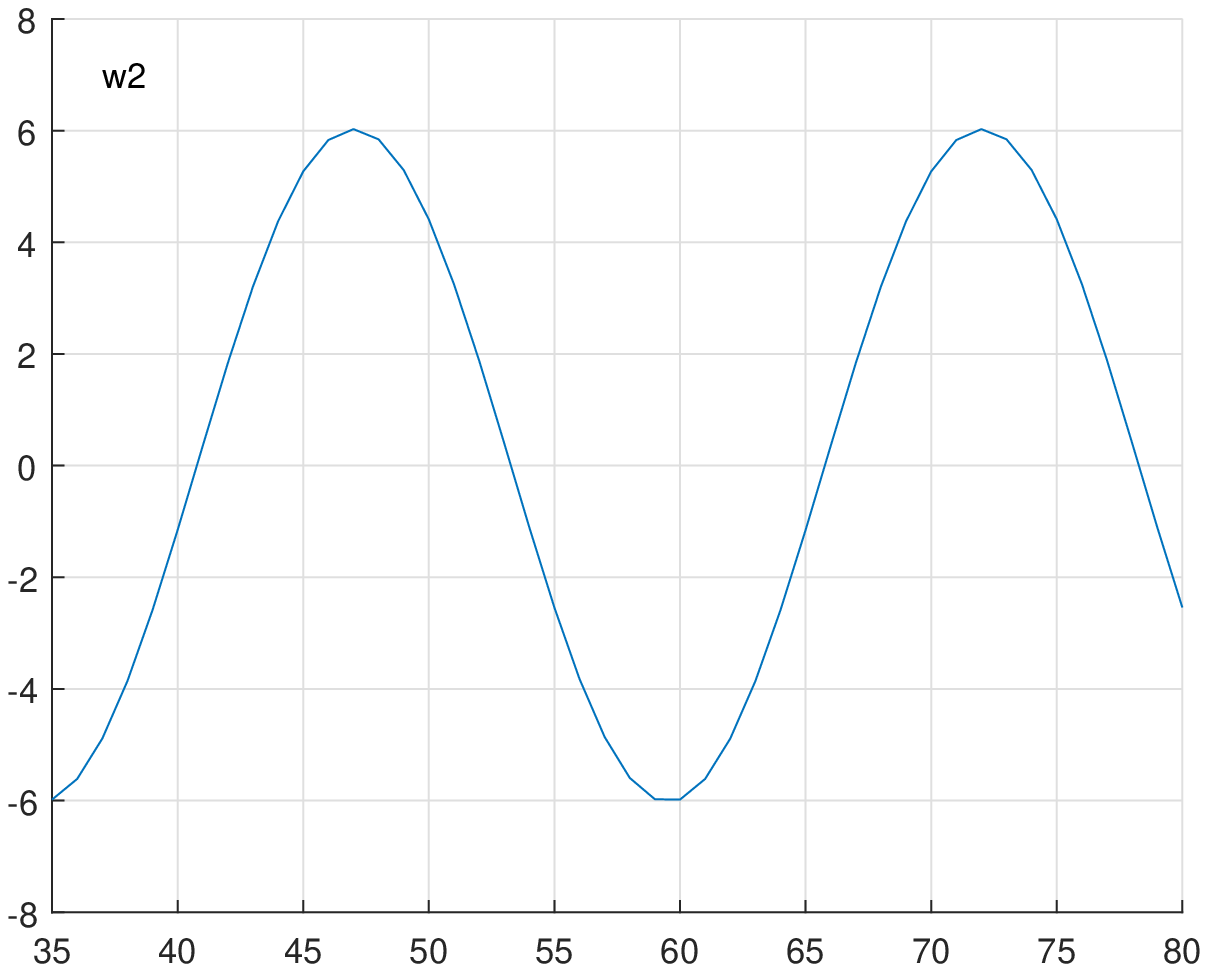}
\label{f2c}\includegraphics[width=0.5\textwidth]{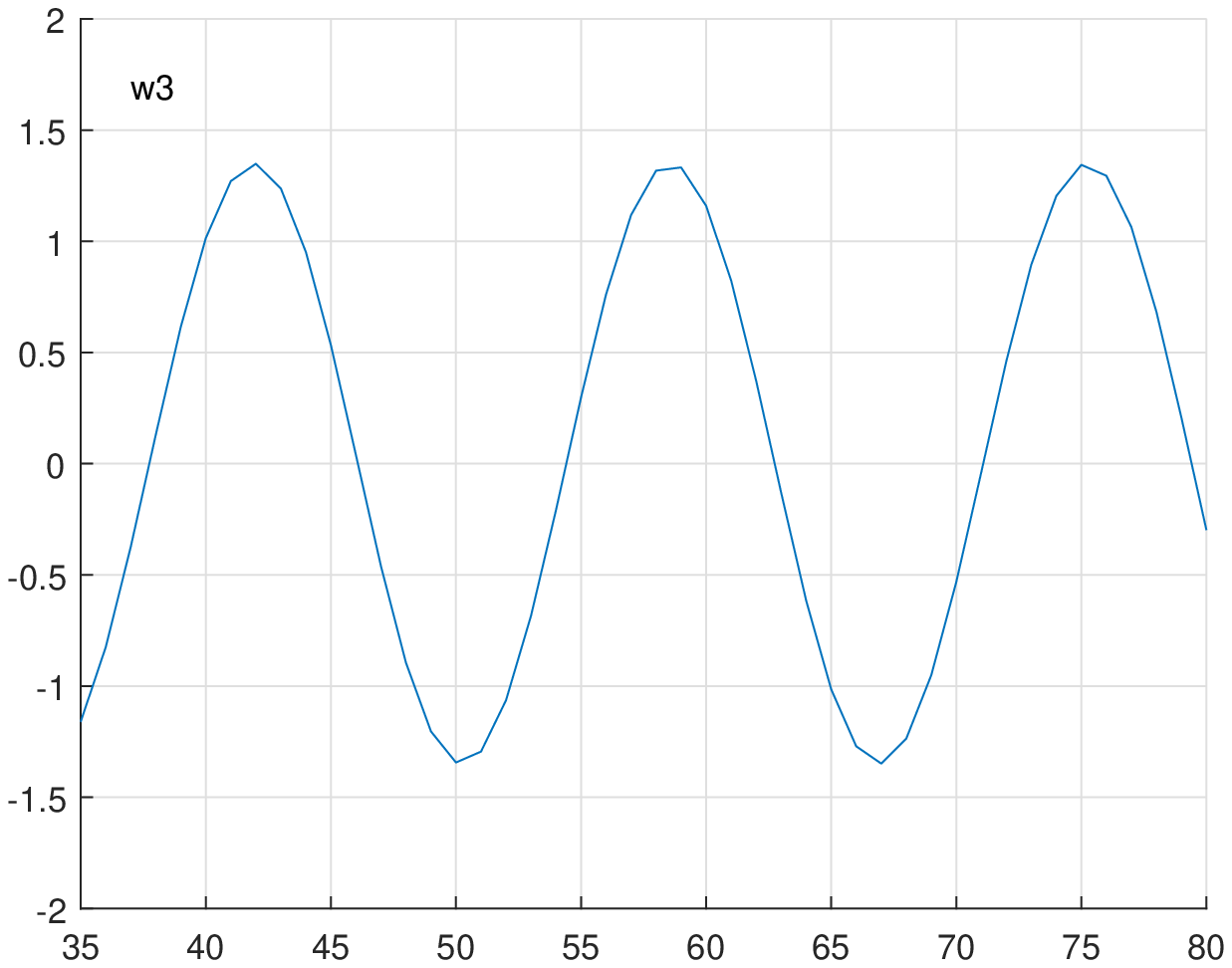}
\label{f2d}\includegraphics[width=0.5\textwidth]{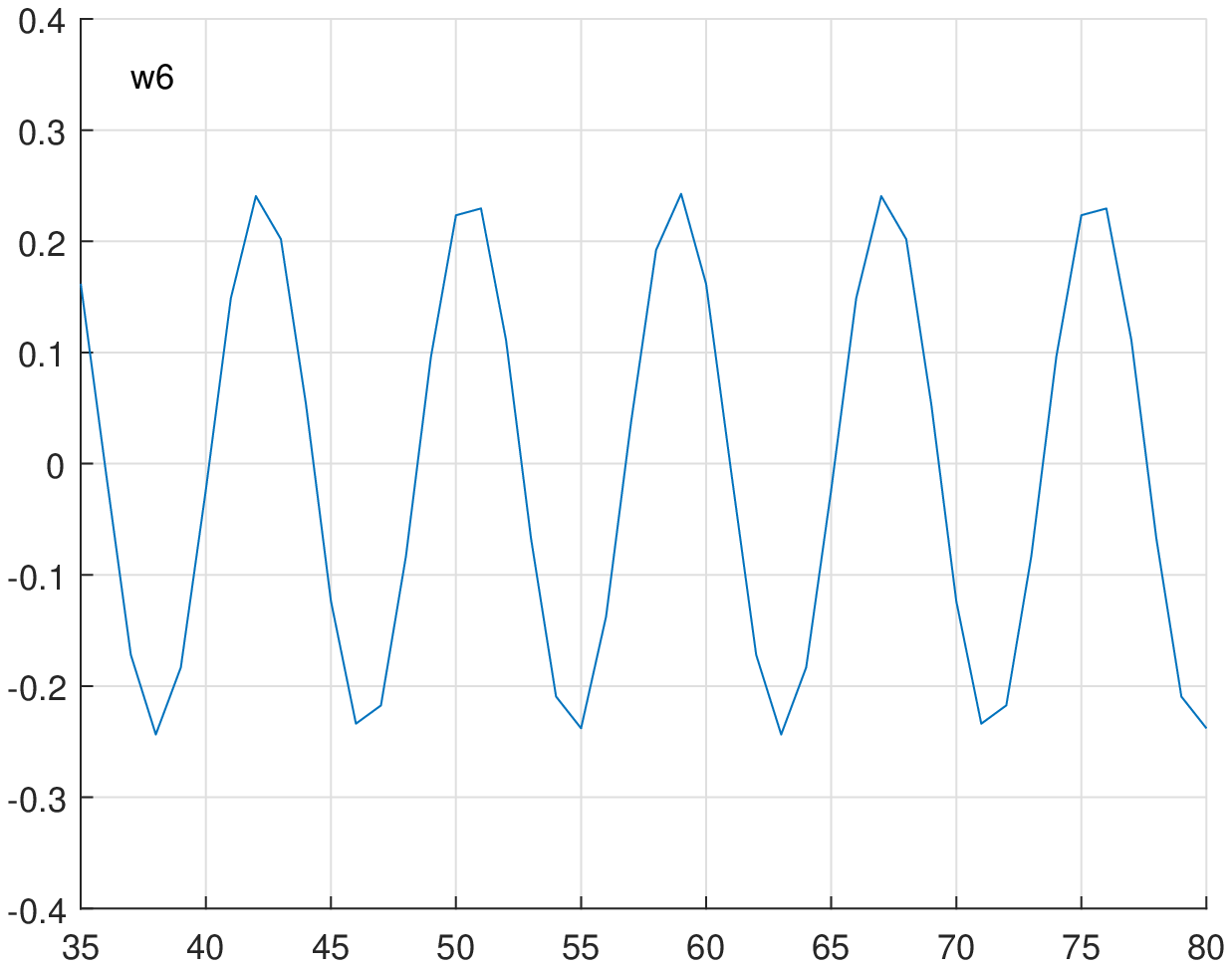}
\caption{Modelling the harmonic waves for the economy of France.}
\label{fig:2}   
\end{figure}

Analysis of modelling curves shows that they give 
a qualitative approximation of statistical data. 
Therefore, the proposed method can be used for effective 
forecasting in real macroeconomic systems. Investigating 
the economy of France, we have received 
the following modelling trajectory:
\begin{equation*}
\begin{split}
x(t)&= -7.10705 +3.3429 t+7.2722 \cos(\omega_1 t)-8.898 \sin(\omega_1 t)\\
&\quad +4.4042 \cos(\omega_2 t) -4.1148 \sin(\omega_2 t)-1.3436 \cos(\omega_3 t)\\
&\quad -0.1239 \sin(\omega_3 t)+0.2233 \sin(\omega_6 t)-0.0976 \sin(\omega_6 t)\;. 
\end{split}
\end{equation*}
Let us analyse the harmonics that are present in this modelling trajectory and, 
therefore, have an influence on the macroeconomic system of France.
Figure~\ref{fig:2} shows the graphics of the first, second, third, and 
sixth harmonics within the interval $t \in [35, 80]$, that is, from years
2000 to 2045.
The analysis shows that since $2000$,
the first three harmonics are in the lifting phase. 
Since $2003$, it is supplemented by the $6$th harmonic lifting phase. 
So this period can be considered as the ``golden era'' of France.
But in $2008$, the $1$st, $3$rd, and $4$th harmonics entered a phase of decline. 
As a result, at the end of $2008$, the economic crisis has begun. 
The cause of the crisis was the fact that all the basic harmonics 
entered into the phase of decline. It is supplemented by the $2$nd 
harmonic since $2012$. Further analysis shows the $6$th harmonic 
will enter the growth phase in $2020$. The $3$rd harmonic is in 
lifting phase now and until $2024$, when the wave starts to decline.
The analysis of the $2$nd wave shows that it is now in a phase of decline.
Changes are expected to begin in $2025$, when the $2$nd wave enters 
the growing phase until $2040$. However, significant changes 
should be expected starting with the year $2034$, 
when the Kondratieff wave will enter the phase of lifting.
In $2034$, it will be supplemented with the $3$rd harmonic, 
in $2037$ with the $6$th harmonics. Therefore, after that time, 
we observe a gradual growth of the economic development of France.


\section{Conclusions}
\label{sec:5}

The constructed mechanism, for simulation and forecasting dynamics 
of a macroeconomic system, allows to establish interconnections 
between individual economy sectors. The identification algorithm 
of the structural model of the inter-branch balance can be used 
for efficient allocation of resources in the formation of 
relationships between different branches. The resulting trajectories 
of issues and non-productive consumption have high imitation 
and forecast properties. The developed method can also be used 
for the analysis and forecasting development of real macroeconomic systems.


\begin{acknowledgement}
The approach used in this paper is based on previous work 
carried out by O.~Kostylenko under the supervision 
of O.~M.~Nazarenko \cite{Olena01,Olena02}, 
which was awarded a diploma at ``All-Ukrainian competition 
of students' research papers''. The research was supported 
by the Portuguese national funding agency for science, 
research and technology (FCT), within the 
Center for Research and Development in Mathematics 
and Applications (CIDMA), project UID/MAT/04106/2019.
Kostylenko is also supported by the Ph.D. 
fellowship PD/BD/114188/2016.
\end{acknowledgement}



\end{document}